\documentclass[12pt,a4paper]{amsart}
\usepackage{amssymb,amsmath,amsthm}
\usepackage{graphicx}
\usepackage{epsfig}
\usepackage{tikz}

\long\def\onefigure#1#2{
\begin{figure*}[tbp]
\begin{center}
#1
\end{center}
\caption{#2}
\end{figure*}
} 

\newcommand{\lipefig}[2]  
{\onefigure{\mbox{\psfig{file=#1.eps}}}{\label{f:#1} #2} }

\newtheorem{theorem}{Theorem}[section]
\newtheorem{lemma}{Lemma}[section]
\newtheorem{corollary}{Corollary}[section]

\newtheorem{claim}{Claim}[section]

\newcommand{\de}{\delta}
\newcommand{\al}{\alpha}
\newcommand{\be}{\beta}
\newcommand{\ga}{\gamma}

\newcommand{\La}{\Lambda}

\newcommand{\eps}{\varepsilon}
\newcommand{\Z}{\mathbb{Z}}
\newcommand{\zz}{{\mathbb{Z}^d}}
\newcommand{\rr}{\mathbb{R}^d}
\newcommand{\R}{\mathbb{R}}

\newcommand{\inn}{\mathrm{int\;}}

\newcommand{\vol}{\mathrm{vol\;}}
\newcommand{\vold}{\mathrm{vol_{d-1}\;}}

\numberwithin{equation}{section}

\begin{document}

\title{Cells in the box and a hyperplane}
\author{Imre B\'ar\'any and P\'eter Frankl}

\begin{abstract} It is well known that a line can intersect at most $2n-1$ cells of the $n \times n$ chessboard. Here we consider the high dimensional version: how many cells of the $d$-dimensional $n\times \ldots \times n$ box can a hyperplane intersect? We also prove the lattice analogue of the following well-known fact. If $K,L$ are convex bodies in $\R^d$ and $K\subset L$, then the surface area of $K$ is smaller than that of $L$.
\end{abstract}

\maketitle

\section{Introduction and main result}\label{sec:introd}

It is well-known that a line can intersect the interior of at most $2n-1$ cells of the $n \times n$ chessboard. What happens in high dimensions? This is the question we address here.

Write $Q_n=Q_n^d=[0,n]^d$, $Q^d=Q_1^d$ so $Q_n^d=nQ^d$. Let $e_1,\ldots,e_d$ be the standard basis vectors of $\rr$ and $\zz$. For $z =(z_1,\ldots,z_d) \in\zz$ define the unit cube
\[
C(z)=\{x=(x_1,\ldots,x_d) \in \rr: z_i \le x_i\le z_i+1, i\in [d]\}
\]
that we are going to call a $\sl cell$ in this paper. Here $[d]$ stands for the set $\{1,2,\ldots,d\}$. For $v \in \rr$, ($v\ne 0$) let $A(v,t)$ denote the hyperplane $\{x\in \rr: vx=t\}$ where $vx$ is the scalar product of the two vectors. Define $N^d(n)$ as the maximal number of cells in $Q_n^d$ that a hyperplane $A(v,t)$ can intersect properly, meaning that $A(v,t)\cap \inn C(z)\ne \emptyset$.

It is well-known that $N^2(n)=2n-1$. Variants of this result have appeared as olympiad problems in several countries. J\'ozsef Beck used a slightly stronger version of this fact to answer questions of Dirac, Motzkin, and Erd\H os in a seminal paper \cite{Beck}. In a companion paper \cite{BarFra} we show that $N^3(n)=\frac 94 n^2+O(n)$. Here we determine the asymptotic behaviour of $N^d(n)$.

We need some definitions. We let $|v|$ resp. $|v|_1$ denote the $\ell_2$ and $\ell_1$ norm of the vector $v\in \rr$. Set
\[
V_d(v)=\frac { |v|_1}{ |v|}\max_{t\in \R} \vold (A(v,t)\cap Q^d),
\]
and
\[
V_d= \max\{ V_d(v): v\in \rr, v\ne 0, t\in \R\}.
\]
It is a consequence of the Brunn-Minkowski theorem (cf \cite{Schn} and the proof of Lemma~\ref{l:slicevol} below) that for fixed $v$ the quantity $\vold (A(v,t)\cap Q^d)$ is maximal when $A(v,t)\cap Q^d$ is the central section of $Q^d$, that is  $A(v,t)$ contains the centre of $Q^d$ which is the point $e/2$ where $e=e_1+\ldots +e_d$. In this case of course $t= ev/2$. It is known that
\[
1\le \vold (A(v, ev/2)\cap Q^d\le \sqrt 2,
\]
the upper bound is a famous result of Keith Ball \cite{Ball}, the lower bound is trivial. This implies that
\[
\sqrt d \le V_d \le \sqrt{2d}.
\]
It is known (see \cite{Ali} or \cite{Ball}) that the sequence $V_2,V_3,\ldots $ is increasing, $V_2=2, V_3=\frac 94, V_4=\frac 83$ etc and its limit is $\sqrt{\frac {6d}{\pi}}$. We conjectured that the vector $v=e$ gives the maximum in the definition of $V_d$. This has been recently proved by Iskander Aliev~\cite{Ali}. Our main result is

\begin{theorem}\label{th:main} $N^d(n)=V_dn^{d-1}(1+o(1))$.
\end{theorem}

In Section~\ref{sec:outline} we give an outline of the proof.

\smallskip
From now on we assume that $v\in \rr$ is a unit vector, i.e., $|v|=1$, and $v\ge 0$, the latter goes without loss generality because of symmetry. Define the (open) strip
\[
S(v,t)=\{x\in \rr: t-ev <vx < t\}.
\]
Clearly
\[
N^d(n)=\max_{v,t}|S(v,t)\cap Q_n^d\cap\zz|.
\]
So we are to determine the number of lattice points in the convex set $S(v,t)\cap Q_n^d$. But this convex set is very thin in one direction (in the direction of $v$) and standard methods do not seem to work. In Section~\ref{sec:bdry} we introduce a novel approach to deal with such cases.

\medskip
 Our result extends to any convex body (convex compact set with non-empty interior) $K \subset \rr$. We define $V(K)=\max \{|v|_1 \vold (K\cap A(v,t)): v\in \rr, |v|=1, t\in \R\}$ and consider the lattice $\frac 1n\zz$.
Write $N(K,n)$ for the maximal number of cells contained in $K$ that a hyperplane can intersect properly (in the same sense as earlier). A cell in this case is $\frac1n C(z)$ with $z \in \zz$. With this notation $N^d(n)=N(Q^d,n)$. Theorem~\ref{th:main} extends to this case as follows.

\begin{theorem}\label{th:mainK} $N(K,n)=V(K)n^{d-1}(1+o(1))$.
\end{theorem}

The proof goes along the same lines as that of Theorem~\ref{th:main} and is therefore omitted.

\bigskip
\section{Inside cells and boundary cells}\label{sec:bdry}

For a general convex body $K$ in $\rr$ a metatheorem says that $\vol K$ is approximately equal to $|K \cap \zz|$, that is
\[
\vol K \approx |K\cap\zz|,
\]
valid when $K$ is well positioned with respect to $\zz$. But this is not necessarily the case with $S(v,t)\cap Q_n^d$. We are going  to well-position it or rather choose a suitable basis of $\zz$ in which $S(v,t)\cap Q_n^d$ is well positioned. We start out more generally.

Let $K \subset \rr$ be a convex body. A cell $C(z)$, $z \in \zz$ called {\sl inside} if $C(z)\subset K$, {\sl outside} if $C(z)\cap K=\emptyset$, and {\sl boundary} otherwise. The following result is going to be useful in other cases as well. It is similar to the well-known fact that the surface area of a convex subset of a convex set $K$ is smaller that the surface area of $K$ itself. To our surprise we couldn't find it anywhere in the literature.\footnote{As we have learned recently, this result was proved by Marek Lassak in~\cite{Lassak} in 1988.  In fact, Lassak's result is more general than ours and is used for a different purpose. His proof, just like ours, goes by a homotopy argument.}

\begin{theorem}\label{th:bdrycell} Assume $K,L$ are convex bodies in $\rr$ and $K\subset L$. Then
\[
| \rm{boundary\; cells\; of\; }{\it K} | \le  | \rm{boundary\; cells \;of \;}{\it L}| .
\]
\end{theorem}

We prove this theorem in Section \ref{sec:bdryproof}.

Now we return to the generic convex body $K$. Since $K$ contains all inside cells and is contained in the union of inside and boundary cells, we have
\[
| \mbox{inside cells of }K|  \le \vol K \le  | \mbox{inside or boundary cells of }K |.
\]
It is not hard to check that
\[
| \mbox{inside cells of }K| \le  |K \cap \zz |  \le  | \mbox{inside or boundary cells of }K|,
\]
implying that
\begin{equation}\label{eq:elem}
\left| \vol K -|K\cap \zz| \right| \le | \mbox{boundary cells of }K |.
\end{equation}

Given a basis $F=\{f_1,\ldots,f_d\}$ of $\zz$ we define the $F$-{\sl box} with parameters $\al,\be \in \rr$ as
\[
B(\al,\be,F)=\{x =\sum_1^d x_if_i \in \rr: \al_i\le x_i \le \be_i, i\in [d]\}.
\]
This is a parallelotope. We of course assume that $\al_i \le \be_i$ for all $i$. The minimal box containing $K$ is denoted by $B(K,F)$, this is the $F$-box  $B(\al,\be,F)$ with all $\al_i$ maximal and $\be_i$ minimal under the condition that $K \subset B(\al,\be,F)$. We will make use of the following theorem of B\'ar\'any and Vershik from \cite{BarVer}, cf \cite{TaoVu} as well.

\begin{theorem}\label{th:BarVer} For every convex body $K$ in $\rr$ there is a basis $F$ such that
\[
\vol B(K,F) \ll_d \vol K.
\]
\end{theorem}

The notation $\ll_d$ means, as usual, that the quantity on the LHS is smaller than the one on the RHS times a positive constant that only depends on $d$. When it is clear from the context, we will use $\ll$ instead of $\ll_d$.  Of course one can use $F$-cells (i.e. basic parallelotopes in the basis $F$) and call them inside, outside, and boundary $F$-cells with respect to $K$. Then inequality (\ref{eq:elem}) becomes
\begin{equation}\label{eq:Felem}
\left| \vol K -|K\cap \zz| \right| \le  | \mbox{boundary }F\mbox{-cells of }K |.
\end{equation}

This inequality extends to any lattice $\La$ and a basis $F$ of $\La$ in the following form:
\begin{equation}\label{eq:FelemLa}
\left| \frac 1 {\det \La}\vol K -|K\cap \La| \right| \le  | \mbox{boundary }F\mbox{-cells of }K |.
\end{equation}

We need a non-degeneracy condition on $K$:
\begin{equation}\label{eq:nondeg}
K\cap \zz \mbox{ contains }d+1 \mbox{ affinely independent vectors}.
\end{equation}
Under this condition and with minimal box $B(K,F)=B(\al,\be,F)$ we have $\al_i \le \lceil \al_i \rceil < \lfloor \be_i \rfloor \le \be_i$ for all $i \in [d]$. Setting $\ga_i=\be_i -\al_i$, $\vol B(K,F)=\prod_1^d\ga_i$. The number of boundary cells of $B(K,F)$ is easy to estimate: it is at most
\[
2\sum_{i=1}^d\prod_{j \ne i}(\ga_j+2)\ll \sum_{i=1}^d\prod_{j \ne i}\ga_j=\vol B(K,F)\left(\frac 1{\ga_1}+\ldots \frac 1{\ga_d}\right).
\]

Combining the previous theorems we have

\begin{theorem}\label{th:basic} Let $K$ be a convex body in $\rr$ satisfying (\ref{eq:nondeg}), and let $F$ be the basis from Theorem~\ref{th:BarVer}. Then
\[
\left| \vol K -|K\cap \zz| \right| \ll \vol K \left(\frac 1{\ga_1}+\ldots +\frac 1{\ga_d}\right).
\]
\end{theorem}

The corresponding version for a general lattice $\La$ says the following. Assume $K$ is a convex body,  $\La$ a lattice in $\rr$, and $K$ contains $d+1$ affinely independent points from $\La$, then there is a basis $F$ of $\La$ such that
\begin{equation}\label{eq:basicLa}
\left| \frac 1{\det \La}\vol K -|K\cap \zz| \right| \ll \frac 1{\det \La} \vol K \left(\frac 1{\ga_1}+\ldots +\frac 1{\ga_d}\right).
\end{equation}

Here just as in Theorem~\ref{th:basic} the parameters $\ga_i$ come from the minimal box $B(K,F)$.

\bigskip
\section{Outline of the proof}\label{sec:outline}

In this section we give a sketch of the proof of Theorem~\ref{th:main}. One main ingredient is Theorem~\ref{th:basic}.

The next section establishes some basic properties of $A(v,t)$ and $S(v,t)$. For instance we show that for fixed $v$ $\vol(S(v,t)\cap Q^n)$ is maximal when $S(v,t)$ is the central strip (Lemma~\ref{l:slicevol}). Write  $S^*(v,t)=S(v,t)\cap Q_n$ for the strip that maximizes, for fixed $v$, the number of lattice points in $S(v,t)\cap Q^n$. We also prove the important but not surprising fact (Lemma~\ref{l:ellips}) that the  convex set $S^*(v,t)$ contains an ellipsoid whose half-axes have lengths of order $n$ apart from one
that has length $|v|_1/2$.

The lower bound in Theorem~\ref{th:main} is simpler and is based on estimating $|S^*(v,t)\cap \zz|$ when $v=z/|z|$ with $z \in \zz$ a primitive vector. In this case the points of $S^*(v,t)\cap \zz$ lie on $|z|_1$ consecutive lattice hyperplanes $A(z,k)$ where $k$ is an integer, and $|A(z,k)\cap \zz|$ is estimated using Theorem~\ref{th:basic} in the form (\ref{eq:basicLa}).

For the upper bound in Theorem~\ref{th:main} we fix a maximizer vector $v=v(n)$ and find a basis $F=F^n=\{f_1,\ldots,f_d\}$ of $\zz$ using Theorem~\ref{th:basic}. This basis is more suitable than the standard one. The main difficulty is to bound $\frac 1{\ga_1}+\ldots +\frac 1{\ga_d}$ on the right hand side of the inequality in Theorem~\ref{th:basic}. Here of course $\ga_i=\ga_i(n)$ for all $i \in [d]$. The upper bound is easy when $\ga_i(n) \to \infty$ for all $i\in [d]$. So we assume that  $\ga_i(n)$ is bounded along a subsequence $n'$ for some $i \in [d]$, for $i=1$, say.

Let $G=G^{n'}$ be the corresponding dual basis, and $g_1(n') \in \zz$ be the corresponding dual basis vector. We show next that $g_1(n')$ is also bounded implying that $g_1(n'')=g$ is a constant (primitive) vector along a further subsequence $n''$. This means that the lattice points in $S^*(v,t)$ lie on $\ga$ consecutive lattice hyperplanes orthogonal to $g$. Here $\ga$ is the floor of $\ga_1(n'')$ which we can assume to be a constant since $\ga_1(n'')$ is bounded. It turns further out that $v(n'')$ tends to $g_0=g/|g|$ because the angle $\phi_{n''}$ between these two vectors is $\ll \frac {\ga}{|g|n''}$.

The next step of the argument is 2-dimensional. Let $\Psi=\Psi_n$ denote the orthogonal projection of $\R^d$ to the 2-plane $\Pi$ spanned by $v(n'')$ and $g$. The projection of the lattice points in $S^*(v,t)$ lie on $\ga $ parallel lines $\ell_h$ (that are $\frac 1{|g|}$ apart) see Figure~\ref{fig:phi}. The projected lattice points on the $h$th line belong to a segment $Y_h$ whose length is $|v(n'')|_1/\sin \phi_{n''}$. We show (Claim~\ref{cl:vertical}) that any line orthogonal to $\ell_h$ intersects at most $\ga +1$ segments $Y_h$, and, more importantly, any such line intersects at most $\ga$ segments $Y_h^*$ where $Y_h^*$ is what you get after deleting a short segment (of length $\sqrt {2d}$) from the left end of $Y_h$.

The number of lattice points in $S^*(v,t)$ is the sum of the lattice points in $\Psi^{-1}(Y_h)$ which is close to $\frac 1{|g|}\vold \Psi^{-1}(Y_h)$ which is close to $\frac 1{|g|}\vold \Psi^{-1}(Y_h^*)$. Estimating the sum of these volumes finishes the proof.

\bigskip
\section{Preparations for the proof of Theorem~\ref{th:main}}\label{sec:prep}

In this section we establish some basic properties of the hyperplane $A(v,t)$ and the strip $S(v,t)$ that give the maximal value of $V_d(v)$. We assume again that $v$ is a unit vector, and suppose without loss of generality that $v \ge 0$, that is, $v_i\ge 0$ for all $i\in [d]$. Actually we can assume that $v_i>0$ for each $i$ because the requirement $A(v,t)\cap \inn C(z)\ne \emptyset$ remains valid even if $v_i$ is modified a little.

For simpler notation we write $A^*(v,t)=A(v,t)\cap Q_n$ and $S^*(v,t)=S(v,t)\cap Q_n$. These intersections  of course depend on $n$, but we suppress this dependence as long as it is not needed. The {\sl central} section is $A^*(v,t_0)$ where $t_0=n|v|_1/2$, this is the section that contains $en/2$, the centre of $Q_n$. The {\sl central} strip is $S^*(v,t_2)$ where $t_2=t_0+|v|_1/2$, this is the strip that is centrally symmetric with centre $en/2$. We will write $A^*(v)$ resp. $S^*(v)$ for the corresponding central section and strip.

\begin{lemma}\label{l:slicevol} For a fixed unit vector $v\in \rr$ $\vol S^*(v,t)$ is maximal for the central strip and
\[
\max_{t \in \R}\vol S^*(v,t)=\vol S^*(v,t_2)=V_d(v)n^{d-1}+O(n^{d-2}).
\]
\end{lemma}

{\bf Proof.} We still assume that $v > 0$ and $|v|=1$. By the Brunn-Minkowski theorem (see \cite{Schn}) the function $t \to \vold A^*(v,t)^{1/(d-1)}$, defined for $t \in [0,n|v|_1]$, is concave. It is also symmetric with respect to $t_0=n|v|_1/2$, and equals zero at the endpoints of the interval $t \in [0,n|v|_1]$. So its maximum is taken at $t_0$, implying that $A^*(v)=A^*(v,t_0)$.  The integral formula
\[
\vol S^*(v,t)=\int_{t-|v|_1}^t \vold A^*(v,s)ds
\]
implies that
\begin{eqnarray*}  \max_{t\in \R} \vol S^*(v,t) &\le& |v|_1 \max_{t\in \R}\vold A^*(v,t)\\
 & =& |v|_1\vold A^*(v)=V_d(v)n^{d-1}.
\end{eqnarray*}

The volume of the central strip is
\[
\vol S^*(v,t_2) =\int_{t_1}^{t_2} \vold A^*(v,t)dt=2\int_{t_1}^{t_0} \vold A^*(v,t)dt
\]
where $t_1=t_0-|v|_1/2$. Concavity implies that on the interval $[t_1,t_0]$
\[
\vold A^*(v,t)\ge \vold A^*(v,t_0)\left(\frac t{t_0} \right)^{d-1}.
\]

We estimate next $D:=|v|_1\vold A^*(v,t_0)-\vol S^*(v,t)$ for $t \in [t_1,t_0]$ :
\begin{eqnarray*}
D&=&2\int_{t_1}^{t_0} \left[\vold A^*(v,t_0)-\vold A^*(v,t)\right]dt\\
    &\le& 2\int_{t_1}^{t_0} \vold A^*(v)\left[1-\left(\frac t{t_0}\right)^{d-1}\right]dt\\
    & \le& |v|_1 \vold A^*(v)\left[1-\left(\frac {t_1}{t_0}\right)^{d-1} \right]\\
        &=&    |v|_1 \vold A^*(v)\left[1-\left(1-\frac 1{2n} \right)^{d-1}\right] < |v|_1 \vold A^*(v)\frac d{2n}.
\end{eqnarray*}
This shows that $\max_{t} \vol S^*(v,t)\ge  V_d(v)\left( 1-\frac d{2n}\right)$. \hfill\qed

\bigskip

Here come the properties of $A(v,t)$ and  $S(v,t)$  that we need. Every $A^*(v,t)$ is contained an a $d-1$-dimensional ball of radius $\ll n$ because $Q_n$ is contained in a ball of radius $\sqrt d n/2$. Fix a unit vector $v$. The {\sl maximizer} is the slice $A^*(v,t)$ that properly intersects the maximal number of cells in $Q_n$ among all $A^*(v,s)$, $s \in \R$. The corresponding $S^*(v,t)$ is also a {\sl maximizer}.

\begin{lemma}\label{l:ball} There is a maximizer $A^*(v,t)$ whose inscribed ball has radius $\gg n$.
\end{lemma}

{\bf Proof.} Recall that $e=e_1+\ldots+e_d$ where vectors $e_1,\ldots,e_d$ form the standard basis of $\R^d$. We can assume by symmetry that the hyperplane $A(v,t)$ satisfies $v>0$ and $t \le ve/2$. $A(v,t)$ contains the (unique) point $a_ie_i$ for all $i\in [d]$, and $a_i>0$, of course. We choose the maximizer hyperplane $A(v,t)$ so that $\min \{a_i: i\in[d]\}$ is maximal. We claim that this maximum is at least $n-1$. Assume that, on the contrary $a_1=\min \{a_i: i\in[d]\} < n-1$. If $A(v,t)$ intersects the cell $C(z)\subset Q_n$, then the hyperplane $A(v,t)+e_1$ intersects the cell $C(z)+e_1$ which lies in $Q_n$, so it intersects at least as many cells as $A(v,t)$. It is easy to check that $A(v,t)+e_1$ contains the (unique) point $a_i'e_i$ with $a_i'>a_i$ for all $i \in [d]$, a contradiction.

Then the $d-1$-dimensional ball inscribed into $A^*(v,t)$ has radius at least $n/d$ as one can see easily.\hfill\qed

\smallskip
We now fix this maximizer $A^*(v,t)$ together with $S^*(v,t)$.

\medskip

\begin{lemma}\label{l:ellips} The maximizer $S^*(v,t)$ contains an ellipsoid with all half-axes having length $\gg n$ apart from one whose length is $|v|_1/2$ which is between $1/2$ and $\sqrt d /2$.
\end{lemma}
\medskip

\medskip
{\bf Proof.} The middle section $A^*(v,t-ev/2)$ of $S^*(v,t)$ contains a $d-1$-dimensional ball of radius $\gg n$. This follows from Lemma~\ref{l:ball} for $n$ large. The width of the strip in direction $v$ is $|v|_1$. \hfill\qed

\bigskip
\section{Lattice points in $A^*(z,h)$}\label{sec:zzz}

Given a primitive vector $z \in \zz$ we are going to estimate the number of lattice points in $A^*(z,h)$ where $h \in \Z$. We will need a more general setting so assume $K$ is a convex subset of $A^*(z,h)$ and we will estimate $|K\cap \zz|$. As $A^*(z,h)$ is $d-1$-dimensional, condition (\ref{eq:nondeg}) asks for $d$ affinely independent points in $K\cap \zz$.

\begin{lemma}\label{l:degen} If $K$ does not satisfy the non-degeneracy condition {\rm (\ref{eq:nondeg})}, then $|K\cap \zz| \ll n^{d-2}$.
\end{lemma}

{\bf Proof.} Under the above conditions the lattice points in $K$ lie on a hyperplane in $A(z,t)$, that is a $d-2$-dimensional affine (lattice) subspace. One can project $K$ orthogonally to a facet of $Q^n$ so that distinct lattice points project to distinct (lattice) points. An induction argument on dimension finishes the proof.\hfill\qed

\begin{lemma}\label{l:nondeg} If $K$ satisfies the non-degeneracy condition {\rm (\ref{eq:nondeg})}, then $|K\cap \zz| \ll \frac 1{|z|}\vold K$.
\end{lemma}

{\bf Proof.} We can apply the general lattice version of Theorem~\ref{th:basic}, i.e., (\ref{eq:basicLa}). The lattice now is $\La=A(z,h)\cap \zz$, it is $d-1$-dimensional and its determinant equals $|z|$, the $\ell_2$ norm of $z$. So there is a basis $F=\{f_1,\ldots,f_{d-1}\}$ of $\La$ such that  $\vold B(K,F) \ll \vold K$. Here $B(K,F)$ is the minimal box in $\La$ containing $K$, and so it is of the form $\{x=\sum_1^{d-1}x_if_i: \al_i \le x_i \le \be_i, i\in [d-1]\}$ with suitable $\al_i,\be_i$. Because of the non-degeneracy assumption  $\ga_i:=\be_i-\al_i\ge 1$. Theorem~\ref{th:basic} shows now that
\[
\left| \frac 1{|z|}\vold K -|K\cap \zz| \right| \ll \frac 1{|z|} \vold K \left(\frac 1{\ga_1}+\ldots +\frac 1{\ga_{d-1}}\right).
\]
As each $\ga_i\ge 1$ this implies the statement.\hfill\qed

\medskip
We assume now that $K \subset A^*(z,h)$ contains a $d-1$-dimensional ball of radius $c_1n$  where $c_1>0$ is a constant depending only on $d$. Of course $K$ lies in a $d-1$-dimensional ball  of radius $\sqrt d n/2$ because $Q_n$ lies in the $d$-dimensional ball of the same radius and centre $en/2$.

\begin{lemma}\label{l:slice} Assume further that $K$ contains $d$ affinely independent points from $\zz$. Then
\[
|K \cap \zz|=\frac 1{|z|} \vold K  \left(1+ |z| O\left(\frac 1n\right)\right),
\]
where the constant in the big Oh notation depends only on $d$.
\end{lemma}

{\bf Proof.} We assume $z\ge 0$ because of symmetry. Again there is a basis $F=\{f_1,\ldots,f_{d-1}\}$ of $\La$ such that  $\vold B(K,F) \ll \vold K \ll n^{d-1}$ where $B(K,F)$ is the minimal box in $\La$ containing $K$ which is of the form $\{x=\sum_1^{d-1}x_if_i: \al_i \le x_i \le \be_i, i\in [d-1]\}$ with suitable $\al_i,\be_i$. Set $\ga_i=\be_i-\al_i$ again and note that  $\vold B(K,F)=|z|\prod_1^{d-1}\ga_i$.

\begin{claim} $n\ll \ga_i |f_i| \ll n$ for every $i\in [d-1]$.
\end{claim}

{\bf Proof.} Let $E$ be the largest volume ($d-1$-dimensional) ellipsoid contained in $B(K,F)$ and define $E^*$ as the blown-up copy of $E$ from its centre by the factor $d-1$. Then $B(K,F)$ is contained in  $E^*$ by the well-known Loewner-John theorem. The volume of $E^*$ is $\ll n^{d-1}$ and $E^*$ contains the ball of radius $c_1n$. This implies that each axis of $E^*$ has length $\gg_d n$ which implies in turn that each axis has length $\ll _d n$. Then the diameter of $E^*$ is $\ll n$, and then so is the diameter of $B(K,F)$ as well. Thus every edge of the parallelotope  $B(K,F)$ has length $\ll n$. These edges are of the form $\ga_if_i$ so $\ga_i|f_i|\ll n$ follows.

On the other hand, the parallelotope $B(K,F)$ contains the ball of radius $c_1n$ so its edges have length at least $c_1n$ showing that $n\ll \ga_i|f|_i$.\hfill\qed

\bigskip
We remark that in view of the claim
\begin{eqnarray*}\label{eq:prod}
n^{d-1} &\gg & \prod \ga_i \prod |f_i|=\frac 1{|z|}\vold B(K,F)  \prod |f_i| \\
&\gg & \frac 1{|z|}n^{d-1}  \prod |f_i|,
\end{eqnarray*}
implying $\prod |f_i|\ll |z|$ and then $|f_i| \ll |z|$ as each $|f_i| \ge 1$.

\medskip
As $\zz \cap K$ contains $d$ affinely independent vectors, Theorem~\ref{th:basic}, or rather its lattice version (\ref{eq:basicLa}) applies. We have using  $|f_i| \ll |z|$ that
\begin{eqnarray*}
 &&\left| \frac 1{|z|} \vold K -|K\cap \zz| \right| \ll \frac 1 {|z|}\vold K \left(\frac 1{\ga_1}+\ldots +\frac 1{\ga_{d-1}}\right)\\
& \ll& \frac 1{|z|}\vold K \left(\frac {|f_1|}n+\ldots +\frac {|f_{d-1}|}n\right) \ll  \frac 1{n} \vold K.
\end{eqnarray*}
So we have indeed
\[
|\zz \cap K|= \frac 1 {|z|}\vold K \left(1+ |z| O\left(\frac 1n\right)\right). \hfill\qed
\]

\bigskip

Set $z_0=z/|z|$ and define
\[
M_d(z,n)= \max_t |\zz \cap S^*(z_0,t)|.
\]

The lattice points in a maximizer $S^*(z_0,t)$ (in the sense used in Lemma~\ref{l:ball}) are all contained in $|z|_1$ consecutive lattice hyperplanes of the form $A(z,h)$. Consequently
\begin{equation}\label{eq:strp}
M_d(z,n)=\max_{k\in \Z}\sum_{h=1}^{|z|_1} |\zz \cap A^*(z,k-h)|.
\end{equation}

\begin{theorem}\label{th:lower} For any primitive vector $z \in\zz$  there is $n_0(z)\in \Z$ such that for all $n>n_0(z)$
\[
M_d(z,n)= n^{d-1}V_d(z_0)+O(n^{d-2}),
\]
where the constant in big Oh notation depends only on $d$.
\end{theorem}

{\bf Proof.} We are to use Lemma~\ref{l:slice} with $K=A^*(z,k-h)$. By Lemma~\ref{l:ball} the maximizer $A^*(z,k)$ contains a ball of radius $\gg n$. It also contains $d$ affinely independent lattice points if $n$ is large enough (depending on $z$). The same applies to all  $A^*(z,k-h)$, $h \in [|z|_1]$ because for large $n$ the slice $A^*(z,k-h)$ is very close to $A^*(z,k)$. We can use  Lemma~\ref{l:slice} in (\ref{eq:strp}) to give that
\[
\sum_{h=1}^{|z|_1} |\zz \cap A^*(z,k-h)|= \sum_{h=1}^{|z|_1} \frac 1{|z|}\vold A^*(z,k-h)\left(1+|z| O\left(\frac 1n\right)\right).
\]
As we have seen $\vold A^*(z,k-h)$ is at most the $d-1$ dimensional volume of the central slice $A^*(z)=A^*(z,t_0)$. So the sum of $\vold A^*(z,k-h)$ for $|z|_1$ consecutive slices is at most $|z|_1\vold A^*(z)$. This sum is maximal when the slices are as close to the central slice as possible. This follows from the concavity of the function $t \to \vold A^*(z,t)^{1/(d-1)}$. The sum of these central slices is estimated the same way as in the proof of Lemma~\ref{l:slicevol}. We omit the details.\hfill\qed

\bigskip
\begin{corollary}\label{cor:lower}
$N^d(n)\ge V_dn^{d-1}(1+o(1))$.
\end{corollary}

{\bf Proof.} Denote by $A^0(v)$ the central section $A(v,t)\cap Q^d$. The function $v\to |v|_1\vold A^0(v)$ (for unit vectors in $\rr$) is continuous. So for any $\eps>0$ we can choose a primitive vector $z\in \zz$ such that $V_d(z_0)\ge V_d-\eps/2$ where $z_0=z/|z|$. Then for all large enough $n$ $M_d(z,n)\ge n^{d-1}V_d(z_0)+O(n^{d-2})\ge n^{d-1}(V_d-\eps/2) +O(n^{d-2}) \ge  n^{d-1}(V_d-\eps)$. \hfill\qed

\bigskip
\section{Proof of the upper bound in Theorem~\ref{th:main}}\label{sec:upper}

Let $S_n=S^*(v,t)$ be the maximizer for $N_d(n)$, of course $v=v(n)$ and $t=t(n)$ but we suppress this dependence as long as possible. We are to show that for every $\eps>0$
\begin{equation}\label{eq:target}
|S_n\cap \zz| \le (V_d+\eps)n^{d-1}
\end{equation}
for all large enough $n$. Fix $\eps>0$.

We claim first that $S_n$ satisfies the non-degeneracy condition ({\ref{eq:nondeg}). Otherwise $S_n \cap \zz$ is contained in a hyperplane of normal $w$ with $we_i\ne 0$ for some $i\in [d]$, $i=d$ say. Projecting the points of $S_n \cap \zz$ to the hyperplane $x_d=0$ we get lattice points on a facet of $Q_n$, and distinct points project to distinct points. No facet contains more than $(n+1)^{d-1}$ lattice points, so $|S_n\cap \zz|\le (n+1)^{d-1}$ which is smaller than $M_d(n)\ge \sqrt d n^{d-1}+O(n^{d-2})$. The last inequality follows from Corollary~\ref{cor:lower} and from $V_d\ge \sqrt d$.

Now Theorem~\ref{th:basic} gives that
\begin{equation}\label{eq:basic}
\left|\vol S_n - |S_n\cap\zz|\right| \ll \vol S_n\left(\frac 1{\ga_1}+\ldots + \frac 1{\ga_{d-1}}\right).
\end{equation}
Here of course $\al_i=\al_i(n), \be_i=\be_i(n)$ and $\ga_i = \ga_i(n)=\be_i(n)-\al_i(n)$.
A simple case is when there is a subsequence $n'$ of the positive integers such
 that $\lim \ga_i(n')=\infty$ for every $i\in [d]$. For simpler writing we use $n$ instead of $n'$.
Then (\ref{eq:basic}) implies that
\[
|S_n\cap\zz|=\vol S_n(1+o(1))\le V_dn^{d-1}(1+o(1)),
\]
so (\ref{eq:target}) holds true indeed.

Assume next that there is a subsequence $n'$ of the previous subsequence such that $\ga_i(n')$ is bounded for some $i \in [d]$, $i=1$ say. We write again $n$ instead of $n'$. Let $G^n=\{g_1^n,\ldots,g_d^n\}$ be the dual basis of $F=F^n$. Set
\[
\al(n)=\min \{g_1^nx: x \in S_n\} \mbox{ and } \be(n)=\max \{g_1^nx: x \in S_n\}.
\]
Of course $\be(n)-\al(n)= \ga_1(n)$ and $\ga_1(n)$ is bounded. So along another subsequence (to be denoted invariably by $n$) $\lim \be(n)-\al(n)= \ga$ for some $\ga \ge 0$.

We claim now that the corresponding dual basis vector $g_1^n$ is also bounded. This is simple again: otherwise the width of $S_n$ in direction $g_1^n$ is $\ga/|g_1^n|$ which tends to zero as $n$ goes to infinity. But $S_n$ contains a ball of radius $\gg 1$ (by Lemma~\ref{l:ellips}), a contradiction. This implies that along a further subsequence $g_1^n$ is equal to a fixed primitive vector, $g$, say.

Define the strip
\[
T_n=\{x\in \rr: \al(n)\le gx \le \be(n)\}, \mbox{ then }S_n\cap\zz  \subset T_n
\]
because of the definition of $\al(n)$ and $\be(n)$. Set $g_0=g/|g|$. Let $\phi_n$ be the angle between $g$ and $v(n)$, so $\cos \phi_n=v(n)g_0$. Define $\Psi: \rr \to \Pi_n$ as the orthogonal projection to the 2-dimensional plane spanned by $v(n)$ and $g$.
Note that here we can assume $g\ne v(n)$ since a minute change of $v(n)$ does not influence what cells the hyperplane $A(v(n),t)$ intersects.

\begin{claim} Along the present subsequence $\phi_n\ll \frac {\ga}{|g|n}$ and so $v(n) \to g_0$.
\end{claim}

\begin{figure}
\centering
\includegraphics[scale=0.9]{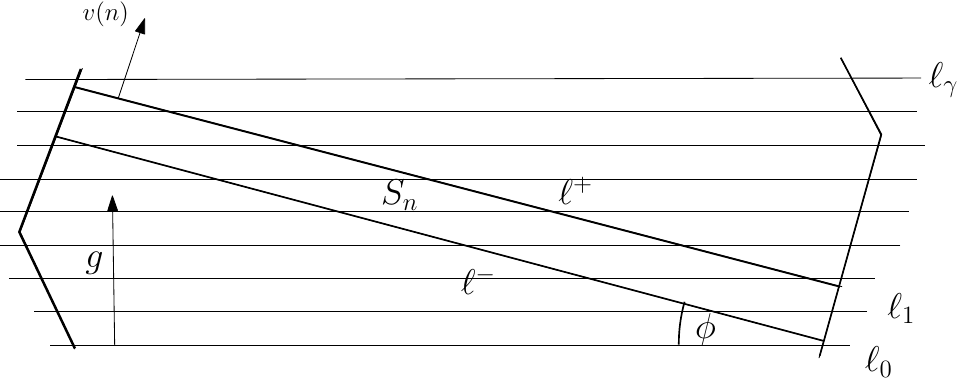}
\caption{$v(n)$ tends to $g_0$}
\label{fig:phi}
\end{figure}

{\bf Proof.}  We drop the subscript $n$ if possible. $\Psi(Q_n)$ is a centrally symmetric convex polygon. The $\Psi$-image of the lattice hyperplane $A(g,\lceil \al(n)\rceil+h)$ is the line $\ell_h$ on $\Pi_n$, represented by a horizontal line on Figure~\ref{fig:phi}, $h=0,1,\ldots,\ga$. Here we take the upper integer part of $\al(n)$ because we need lattice hyperplanes.
We should also take $h=0,1,\ldots,\lfloor \ga \rfloor$ because $\ga$ may not be integer. But for simpler writing we keep it $\ga$ now and in what follows.

The $\Psi$-image of the two hyperplanes bounding $S_n=S(v(n),t(n))$ are the lines $\ell^+$ and $\ell^-$ on Figure~\ref{fig:cov}. Their distance is $|v|_1$. The length of the segments $\ell^+\cap \Psi (Q^n)$ and $\ell^+\cap \Psi (Q^n)$ is $\gg n$ because $S_n$ contains the ellipsoid from Lemma~\ref{l:ellips} and $S_n\subset T_n$. So with $\phi=\phi_n$
\[
\sin \phi \ll \frac {\ga}{n|g|}.
\]\hfill\qed

\smallskip
Define now for $h=0,1,\ldots,\ga$ the $d-1$-dimensional convex polytope
\[
P_h^n=S_n\cap A(g,\lceil \al(n)\rceil +h),
\]
Every lattice point in $S_n$ belongs to some $P_h^n$. The proof of the upper bound on $M_d(n)$ is based on estimating $\sum_0^{\ga}|P_h^n\cap \zz|$.

Define the map $\Phi=\Phi_n: \rr \to \rr$ via $\Phi(x)=x/n$. Then $\Phi(P_h^n)$ is convex compact sets in $Q^d$ for all $h\in \{0,1,\ldots.\ga\}$. We use the Blaschke selection theorem, see for instance~\cite{Schn}:  along a subsequence (denoted by $n$ again) $\Phi(P_h^n)$ tends to a convex polytope $P_h$  for $h\in \{0,1,\ldots,\ga\}$. Note also that each $P_h$ lies in $A(g,t)\cap Q^d$ for some fixed $t$.

\smallskip
Let $I$ denote the set of $h\in \{0,\ldots,\ga\}$ with $\vold \Psi(P_h)>C_0$ where $C_0>0$ will be specified later. Write $J_1$ resp. $J_0$ for those $h \notin I$ for which $P_h^n$ does (does not) contain $d$ affinely independent vectors from $\zz$. We are going to estimate $|P_h^n\cap \zz|$ separately for $h$ in $I$ and $J_0$ and $J_1$.

\smallskip
When $h \in J_0$ Lemma~\ref{l:degen} applies and gives $|P_h^n \cap \zz|\ll n^{d-2}$. The total contribution of such $P_h^n$s to $|S_n\cap \zz|$ is at most $\ll |J_0|n^{d-2}$.

\smallskip
For $h \in J_1$  Lemma~\ref{l:nondeg} shows that
$|P_h^n\cap \zz| \le C_d \frac 1{|g|} \vold P_h^n$. Here $C_d>0$ is the constant implicit in the $\ll$ notation.  The total contribution of such $P_h^n$s to $|S_n\cap \zz|$ is at most $\ll |J_1| \frac {C_dC_0}{|g|}n^{d-1}\le |J_1|C_dC_0n^{d-1}$.

\smallskip
For $h \in I$ let $E_h^n$ be the ellipsoid of largest volume inscribed in $P_h^n$ with half-axis of length $a_1,\ldots,a_{d-1}$. The Loewner-John theorem implies that
\[
\vold E_h^n \ge (d-1)^{-(d-1)}\vold P_h^n \ge C_0 \left(\frac n{d-1}\right)^{d-1}.
\]
Also $\vold E_h^n=\kappa_{d-1}\prod_1^{d-1}a_i$ where $\kappa_{d-1}$ is the volume of the $d-1$-dimensional unit ball. As each $a_i\le \sqrt dn$, the minimal $a_i\gg C_0n$. So $P_h^n$ contains a ball of radius $\gg C_0n$. It is also clear that for large enough $n$ $P_h^n$ contains $d$ affinely independent points from $\zz$. So we can apply Lemma~\ref{l:slice}: for $h \in I$
\[
|P_h^n \cap \zz| \le \frac 1{|g|}\vold P_h^n\left(1+ |g|O\left(\frac 1n\right)\right),
\]
showing that the total contribution of these $P_h^n$s to $|S_n\cap \zz|$ is at most
\[
\frac 1{|g|}\sum_{i\in I} \vold P_h^n\left(1+ |g| O\left(\frac 1n\right)\right).
\]

\begin{lemma}\label{l:aux} With the previous notation
\[
\frac 1{|g|}\sum_0^{\ga} \vold P_h ^n\le V_d(g)n^{d-1}\left(1+o(1)\right).
\]
\end{lemma}

We postpone the proof to the next section. We show now how to complete the proof of Theorem~\ref{th:main} using this lemma.

The number of lattice points in $S_n$ is  $V_d(g)n^{d-1}(1+o(1)) \le V_dn^{d-1}+\frac 12 \eps n^{d-1}$ if $n$ is large enough plus an error term of the form
\[
|J_0|n^{d-2}+|J_1|C_dC_0n^{d-1}
\]
times a constant depending only on $d$. Here $|J_1|,|J_0|\le \ga$, and $g$ and $\ga$ are fixed. So choosing $C_0>0$ small enough the error term becomes smaller than $\frac 12 \eps n^{d-1}$.\hfill\qed

\bigskip
\section{Proof of Lemma~\ref{l:aux}}

\begin{figure}
\centering
\includegraphics[scale=0.8]{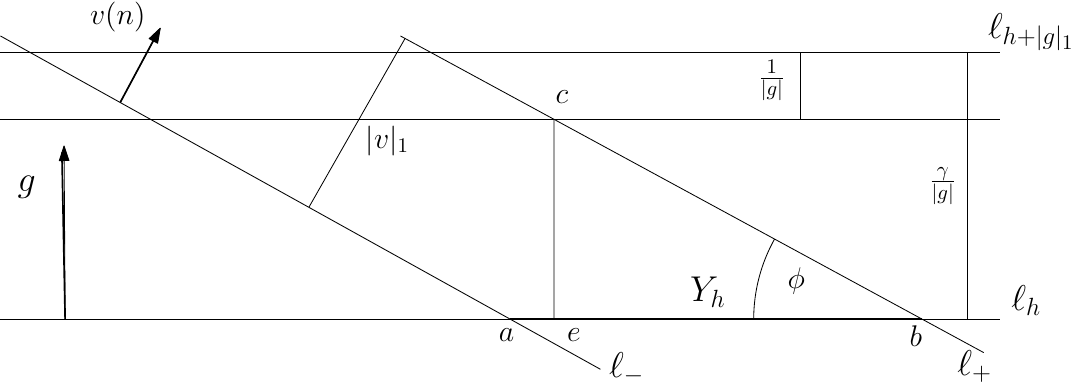}
\caption{Projection onto $\Pi$}
\label{fig:cov}
\end{figure}

We first note that $(v(n)-g_0)^2=2-2\cos \phi=2\sin^2 \phi/(1+\cos \phi)$. Set $Y_h=\ell_h\cap \Psi(S_n)$; this is a segment whose length is $|v|_1/\sin \phi$. Let $Y_h^* \subset Y_h$ be the segment that you get after deleting the segment of length $\sqrt{2d}$ from the left end of $Y_h$.

\begin{claim}\label{cl:vertical} A vertical line intersects at most $|g|_1+1$ segments $Y_h$ and at most $|g|_1$ segments $Y_h^*$, $h=0,1,\ldots,\ga$.
\end{claim}

The {\bf proof} is elementary plane geometry using the fact that $v(n)$ and $g_0$ are very close by. We assume $v(n)>0$, then $g\ge 0$ as well and $|v(n)|_1=v(n)e$, $|g|_1=ge$. Assume $\ell^-$ intersects $\ell_h$ in point $a$, and $\ell^+$ intersects $\ell_h$ resp. $\ell_{h+|g|_1}$ in points $b$ and $c$, and let $e$ denote the orthogonal projection of $c$ to $\ell_h$. We consider $a,b,e$ as real numbers on the $x$ axis.  The length of $Y_h$ is $b-a=ve/\sin \phi$, and $b-e=|g_1|/(|g|\tan \phi)=g_0e/\tan \phi$ and
\begin{eqnarray*}
e-a&=&\frac {ve}{\sin \phi}-\frac {g_0e}{\tan \phi}=\frac 1{\sin \phi}\left({ve-g_0e \cos \phi}\right) \\
     &=&\frac 1{\sin \phi} \left[(v-g_0)e+g_0e(1-\cos \phi)\right]\\
     &\le& \frac 1{\sin \phi}\left( \frac {\sqrt 2 \sin \phi}{\sqrt{1+\cos \phi}}\sqrt d+\frac {\sin^2 \phi}{1+\cos \phi}\right)<\sqrt {2d},
\end{eqnarray*}
as one can check easily. This implies that $Y_h^*$ is contained in the interval $[e,b]$. Moreover, a vertical line intersecting the segment $[a,e]$ intersect $Y_h,Y_{h+1},\ldots, Y_{h+|g|_1}$ but no other $Y_i$. And a vertical line intersecting the segment $(e,b]$ intersects $Y_h,\ldots,Y_{h+|g|_1-1}$ but no other $Y_i$. \hfill\qed
\smallskip

The claim implies what we need. Note that $P_h^n=\Psi^{-1}(Y_h)\cap Q^n$, and define $P_h^{n*}=\Psi^{-1}(Y_h^*)\cap Q^n$. Then $P_h^{n*} \subset P_h^n$ and evidently
\[
\vold P_h^n-\vold P_h^{n*} =O(n^{d-2}).
\]
Recalling that $\Phi(x)=x/n$ we have
\[
\sum_0^{\ga}\vold P_h^{n*}=n^{d-1}\sum_0^{\ga}\vold \Phi(P_h^{n*}).
\]

The sets $\Phi(P_h^{n*})$ tend to a set $P^h \subset A(g,t)\cap Q^d$ for the same $t$ as before, so $\vold \Phi(P_h^{n*})=n^{d-1}\vold P^h(1+o(1))$. The sets $P^h$, for $h=0,1,\ldots,\ga$ cover $A(g,t)\cap Q^d$ at most $|g|_1$ times. So their total $d-1$-volume is at most $|g|_1 \vold A(g,t)\cap Q^d$. Thus
\begin{eqnarray*}
\sum_0^{\ga}\vold P_h^n &\le& n^{d-1}\sum_0^{\ga}\vold \Phi(P_h^{n*})+O(n^{d-2}) \\
            & \le &    n^{d-1}\sum_0^{\ga} \vold P^h(1+o(1))\\
            & \le &n^{d-1} |g|_1 \vold (A(g,t)\cap Q^d)(1+o(1)).
\end{eqnarray*}

So indeed
\begin{eqnarray*}
\frac 1{|g|}\sum_0^{\ga}\vold P_h^n &\le& \frac {|g|_1}{|g|} \vold (A(g,t)\cap Q^d)(1+o(1))\\
      &=& V_d(g_0)(1+o(1))
\end{eqnarray*}
because $\frac {|g|_1}{|g|} \vold (A(g,t) \cap Q^d)\le V_d(g_0)$ by the definition of  $V_d(g_0)$.\hfill\qed

\bigskip
\section{Proof of Theorem~\ref{th:bdrycell}}\label{sec:bdryproof}

We construct a homotopy $t\to K_t$ where $t \in [0,1]$, $K_t$ is a convex body in $\rr$ satisfying $K_0=K,\; K_1=L$ and the monotonicity condition $K_t \subset K_s$ for every pair $t<s$. Because of monotonicity boundary cells of $K_t$ may become inside cells for $K_s$, and the point in the argument is that whenever a boundary cell is lost, another one emerges.

The simplest homotopy is $K_t=(1-t)K+tL$, and this works under the following non-degeneracy condition:

\smallskip
(*) whenever $w \in \partial K_t\cap \zz$, then $w \notin K_s$ for $s<t$ and $w \in \inn K_s$ for all $s>t$, and $K_t$ has an outer normal  $u$ at $w \in \partial K_t$ with no coordinate zero.
\smallskip

Under this condition the proof is easy. As $t$ increases, a cell $C(z)$, say, is boundary for $K_t,\; t<t_0$ just slightly smaller than $t_0$ but $C(z) \subset K_{t_0}$ and so it becomes inside for $t>t_0$. Then there is a vertex $w$ of $C(z)$ such that $w \notin K_t$ for $t<t_0$, but of course $w \in K_{t_0}$ and even $w \in \partial K_{t_0}$. Let $H$ be a supporting hyperplane to $K_{t_0}$ at $w$ whose outer normal has no zero coordinate. Then $w \in K_{t_0}$ and $C(z)$ and $K_{t_0}$ are on the same side of $H$. There is a unique cell $C(z')$ (unique because of condition (*)) on the other side of $H$ with $w \in C(z')$. This unique cell was outside for $K_t, \; t<t_0$ and becomes boundary for $K_t$ for $t\in [t_0,t_0+\de)$ for a suitable small $\de>0$. So when the boundary cell $C(z)$ is lost at $t_0$ another boundary cell appears.
Note that $C(z)\cap H=C(z')\cap H=\{w\}$.

We have to check yet that the same cell $C(z')$ can't appear twice. So assume the contrary, that is, there is another cell $C(z^*)$ that is boundary for $K_t$, for $t$ slightly smaller than $t_0$ but $C(z^*)\subset K_{t_0}$ and $C(z^*)$ has a vertex $w^*$ with $w^* \notin K_t$ for $t<t_0$ but $w^*\in \partial K_{t_0}$. We can't have $w=w^*$ here since that would imply $C(z)=C(z^*)$. Then $w$ and $w^*$ are distinct vertices of $C(z')$ and the segment $[w,w^*]$ is on the boundary of both $C(z')$ and $K_{t_0}$. Then $[w,w^*]\cap H=\{w\}$ for the previous hyperplane $H$ supporting $K_{t_0}$ at $w$ with no zero coordinate so $w^*\in K_{t_0}$ can't hold.

\smallskip 
 To guarantee the non-degeneracy condition we proceed first by assuming that $K \subset \inn L$ and that  both $K$ and $L$ have smooth boundaries such that for every unit vector $u$ there is a single point on $\partial K$ resp. on $\partial L$ where the outer normal to $K$ and $L$ is $u$. If this were not the case, we can replace $K,L$ by suitable (and very close to $K$ and $L$) convex bodies satisfying these conditions and having the same inside and boundary cells. With the new $K$ and $L$ the homotopy $K_t=(1-t)K+tL$ has the property that for every unit vector $u$ there is a single point on $\partial K_t$ where the outer normal to $K_t$ is $u$. To see that this is indeed the case, let $x_K$ and $x_L$ be the unique points on the boundary of $K$ and $L$ with outer normal $u$. Then the maximum of $\{ux: x\in K_t\}$ is reached on the unique point $(1-t)x_K+tx_L\in K_t$, and the outer normal to $K_t$ there is $u$.

\smallskip
This condition also guarantees that $K_t$ has no line segment on its boundary. Assume that, on the contrary, $\partial K_t$ contains a line segment and let $u$ be the outer normal to the tangent hyperplane to $K_t$ containing this segment. Then  there is no unique point with outer normal $u$ as every point on the segment has outer normal $u$.

\smallskip
Let's see finally that $K_t$ satisfies condition (*) as well. Assume the cell $C(z)$ is boundary for $K_t$ for
$(t_0-\de,t_0)$ and is inside for $K_{t_0}$. Then there is a vertex $w$ of $C(z)$ on $\partial K_{t_0}$, with outer normal $u=(u_1,\ldots,u_d)$ at $w$ to $K_{t_0}$. Assume some coordinate of $u$ is equal to zero, say $u_1=0$. Either $w+e_1$ or $w-e_1$ is in $C(z)$, say $w+e_1$. Then the segment $[w,w+e_1]$ lies both in $K_{t_0}$ and in $C(z)$, and actually in the boundary of both because the hyperplane $\{x:ux=uw\}$ is tangent to both $K_{t_0}$ and $C(z)$.   \hfill\qed

\bigskip
{\bf Acknowledgements.}  We acknowledge priority of Marek Lassak concerning Theorem~\ref{th:bdrycell}. His result in~\cite{Lassak} is more general than our theorem and was proved much earlier, in 1988.  This piece of research was supported by Hungarian National Research Grants No 131529, 131696, and 133819.

\bigskip

\noindent
Imre B\'ar\'any \\
R\'enyi Institute of Mathematics,\\
13-15 Re\'altanoda Street, Budapest, 1053 Hungary\\
{\tt barany.imre@renyi.hu}\\
and\\
Department of Mathematics\\
University College London\\
Gower Street, London, WC1E 6BT, UK

\medskip
\noindent
P\'eter Frankl \\
R\'enyi Institute of Mathematics,\\
13-15 Re\'altanoda Street, Budapest, 1053 Hungary\\
{\tt peter.frankl@gmail.com}.


\bigskip



\begin{thebibliography}{999}

\bibitem{Ali}
Aliev, I., On the volume of hypercube sections, {\it Acta Math. Hungar.},  \textbf{163} (2021,  547--551.


\bibitem{Ball}
Ball, K. M., Cube slicing in $\R^n$, {\it Proc. Amer. Math. Soc.}, \textbf{97} (1986), 465--473.

\bibitem{BarFra}
B\'ar\'any, I., Frankl, P., How (not) to cut your cheese, {\it Amer. Math. Monthly}. \textbf{128} (2021),  543--552.

\bibitem{BarVer}
B\'ar\'any, I., Vershik, A., On the number of convex lattice
polytopes,  {\it Geom. Functional Analysis}, {\bf 2 } (1992),
381--393.

\bibitem{Beck}
Beck, J., On the lattice property of the plane and some problems of Dirac, Motzkin and
Erd\H os in combinatorial geometry, {\it Combinatorica},  \textbf{3} (1983), 281--297.

\bibitem{Lassak}
Lassak, M., Covering the boundary of a convex set by tiles. {\it Proc. Amer. Math. Soc.},  \textbf{10} (1988), 269--272. 

\bibitem{Schn}
Schneider, R., {\it Convex bodies: the Brunn Minkowski theory}. Second expanded edition. Encyclopedia of Mathematics and its Applications, 151. Cambridge University Press, Cambridge, 2014.

\bibitem{TaoVu}
Tao, T, and Vu, Van H., {\it Additive combinatorics}. Cambridge Studies in Advanced Mathematics, 105. Cambridge University Press, Cambridge, 2010.

\end{thebibliography}
\end{document}